\documentclass{proc-l}
\newtheorem{theorem}{Theorem}[section]
\newtheorem{lemma}[theorem]{Lemma}

\newtheorem{conjecture}[theorem]{Conjecture}

\theoremstyle{definition}
\newtheorem{definition}[theorem]{Definition}
\newtheorem{example}[theorem]{Example}

\theoremstyle{remark}
\newtheorem{remark}[theorem]{Remark}

\numberwithin{equation}{section}

\def\spa{{\mbox{\rm span }\!}}
\def\sp{\spa}

\def\NN{\mathbb{N}}

\def\cH{\mathcal{H}}
\def\cK{\mathcal{K}}

\def\NN{\mathbb{N}}

\def\RR{\mathbb{R}}

\renewcommand{\H}{\mathcal H}

\newcommand{\N}{\mathbb N}

\newcommand{\ip}[2]{\left\langle#1,#2\right\rangle}
\newcommand{\absip}[2]{\left| \left\langle#1,#2\right\rangle \right|}
\newcommand{\norm}[1]{\left\| #1 \right\|}

\begin{document}

\title[Decompositions of frames and
the Feichtinger Conjecture]{A Decomposition Theorem for frames
and the Feichtinger Conjecture}

\author{Peter G. Casazza}
\address{Department of Mathematics,
University of Missouri,
Columbia, MO 65211, USA}
\email{pete@math.missouri.edu}
\thanks{The first author was supported  by NSF Grant DMS 0405376.}

\author{Gitta Kutyniok}
\address{Institute of Mathematics,
Justus--Liebig--University Giessen,
35392 Gies\-sen, Germany}
\email{gitta.kutyniok@math.uni-giessen.de}
\thanks{The second author was supported by DFG Research Fellowship KU 1446/5-1 and by Preis der
Justus-Liebig-Universit\"at Gie{\ss}en 2006.}

\author{Darrin Speegle}
\address{Department of Mathematics and Computer Science,
Saint Louis University,
Saint Louis, MO 63103, USA}
\email{speegled@slu.edu}
\thanks{The third author was supported by NSF Grant DMS 0354957.}

\author{Janet C. Tremain}
\address{Department of Mathematics,
University of Missouri,
Columbia, MO 65211, USA}
\email{janet@math.missouri.edu}

\subjclass{46C05, 42C15, 46L05}

\keywords{Bessel sequence, decomposition, frame, Feichtinger Conjecture, frame sequence,
Kadison-Singer Conjecture, $\omega$-independence, Riesz basic sequence}

\begin{abstract}
In this paper we study the Feichtinger Conjecture in frame theory, which was
recently shown to be equivalent to the 1959 Kadison-Singer Problem
in $C^{*}$-Algebras.
We will show that every bounded Bessel sequence can be
decomposed into two subsets each of which is an arbitrarily
small perturbation of a sequence with a finite orthogonal
decomposition.  This construction is then used to answer two
open problems concerning the Feichtinger Conjecture:
1. The Feichtinger Conjecture is equivalent to the conjecture that every
unit norm Bessel sequence is a finite union of frame sequences.
2. Every unit norm Bessel sequence 
is a finite union of sets each of which is $\omega$-independent for
$\ell_2$-sequences.
\end{abstract}

\maketitle

\section{Introduction}

The Kadison-Singer Problem \cite{KS59} in $C^{*}$-Algebras has remained unsolved
since 1959, thereby defying the best efforts of several of the most talented
mathematicians in our time. Recently, there has been a flurry of activity
around this problem due to a fundamental paper by the first and fourth author
\cite{CT06} (cf. also the longer version joint with M. Fickus and E. Weber \cite{CFTW06}),
which connects the Kadison-Singer Conjecture with many longstanding open
conjectures in a variety of different research areas -- in Hilbert
space theory, Banach space theory, frame theory, harmonic analysis, time-frequency
analysis, and even in engineering -- by proving that these conjecture are
in fact equivalent to the Kadison-Singer Problem.

\smallskip

In this paper we focus on the equivalent version of the Kadison-Singer
Problem in {\em frame theory}, the so-called {\em Feichtinger Conjecture}.
Before elaborating on the history of this
conjecture and the contribution of our paper, let us first recall the basic
definitions and notations in frame theory.

A countable collection of elements
$\{f_i\}_{i \in I}$ is a {\em frame} for a separable Hilbert space $\H$, if
there exist $0 < A \le B < \infty$ (the \emph{lower} and \emph{upper frame bound})
such that for all $g \in \H$,
\begin{equation}
\label{framedef}
A \, \norm{g}^2 \le \sum_{i \in I} \absip{g}{f_i}^2
\le B \, \norm{g}^2.
\end{equation}
A frame $\{f_i\}_{i \in I}$ is {\em bounded}, if $\inf_{i \in I}
\|f_i\|>0$.  Note that $\sup_{i \in I}\|f_i\|<\infty$ follows
automatically by \cite[Proposition 4.6]{Cas00}, and {\em unit norm}, if
 $\|f_i\| = 1$ for all $i \in I$. If $\{f_i\}_{i \in I}$ is a frame
only for its closed linear span, we call it a {\em frame sequence}.
Those sequences which satisfy the upper inequality in (\ref{framedef})
are called {\em Bessel sequences}.
A family $\{f_i\}_{i \in I}$ is a {\em Riesz basic sequence} for $\H$, if
it is a \emph{Riesz basis} for its closed linear span, i.e., if there
exist $0 < A \le B < \infty$ such that for all sequences of scalars $\{c_i\}_{i \in I}$,
\[ A \, \sum_{i \in I} |c_i|^2 \le \left\|\sum_{i \in I} c_i f_i\right\|^2 \le B \,
\sum_{i \in I} |c_i|^2.\]
Finally, a sequence $(f_i)_{i \in I}$ is called {\em $\omega$-independent for $\ell_2$-sequences},
if, whenever $c=(c_i)_{i \in I}$ is an $\ell^2$-sequence of scalars and $\sum_{i\in I}c_{i}f_i= 0$,
it follows that $c = 0$.

\smallskip

Having recalled the necessary definitions and notations, we can now state
the main conjecture we will be addressing in this paper.

\begin{conjecture}[Feichtinger Conjecture]
\label{F_conj}
Every bounded frame can be written as a finite union of Riesz
basic sequences.
\end{conjecture}

Much work has been done on the Feichtinger Conjecture in just the last few
years \cite{G03,CCLV04,CV05,CT06,CFTW06,BS06}. In particular, by employing
the equivalence of the Paving Conjecture to the Kadison-Singer Problem
shown by Anderson in 1979 \cite{And79} and by using the Bourgain-Tzafriri
Conjecture \cite{BT91} which arose from the ``restricted invertibility principle''
by Bourgain and Tzafriri from 1987 \cite{BT87}, the series of papers \cite{CCLV04,CV05,CT06,CFTW06}
proves the equivalence between the Kadison-Singer Problem and the Feichtinger
Conjecture. In \cite{G03} and \cite{BS06} the Feichtinger Conjecture is considered for
special frames such as wavelet and Gabor frames and frames of translates.
In several cases the Feichtinger Conjecture could indeed be verified, thereby
verifing parts of the Kadison-Singer Problem.


Let us now take a closer look at the Feichtinger Conjecture (Conjecture \ref{F_conj}).
It is easily seen, by just normalizing the frame vectors,
 that we may assume in Conjecture \ref{F_conj} that the frame is a unit norm
frame.  It also follows easily that we
only need to assume that the sequence is a bounded Bessel sequence.
That is, by adding an orthonormal basis
to the Bessel sequence we obtain a bounded frame which can be written
as a finite union of Riesz basic sequences if and only if the original
Bessel sequence can be written this way. Thus the Feichtinger Conjecture
reduces to the conjecture that every unit norm Bessel sequence can be
written as a finite union of Riesz basic sequences.

\smallskip

The first main result of our paper concerns a further reduction of the Feichtinger
Conjecture. For this, consider the following
conjecture which intuitively seems to be much weaker than Conjecture \ref{F_conj}.

\begin{conjecture}
\label{our_conj}
Every unit norm Bessel sequence can be written as a finite union of frame
sequences.
\end{conjecture}

However, surprisingly, we will show that both conjectures 
are in fact equivalent.

\begin{theorem}
\label{main1}
The Feichtinger Conjecture is equivalent to Conjecture \ref{our_conj}.
\end{theorem}

Part of the motivation of this result is a result of
Casazza, Christensen, and Kalton \cite{CCK} concerning frames of
translates, which shows that the set of translates of a function in $L^2(\RR)$
with respect to a subset of $\NN$ is a frame sequence if and only if it is a
Riesz basic sequence.

\smallskip

It is generally expected that the Kadison-Singer Problem will turn out to be
false, which calls for positive partial results. Our second main result answers
an open problem concerning weakenings of the famous Paving Conjecture 
of Anderson \cite{And79} which he showed is
equivalent to KS.  
The main question has been whether every
unit norm Bessel sequence is a finite union of $\omega$-independent
sets.  Recall that a set of vectors $\{f_i\}_{i=1}^{\infty}$ 
is $\omega$-independent if $\sum_{i=1}^{\infty}a_if_i = 0$ implies
$a_i=0$ for all $i=1,2,\ldots$.  This concept was defined
in the 1950's by Marc Krein except that he requested that
the above implication hold for sequences $\{a_i\}_{i\in I}\in \ell_2$.
 Since it is known \cite{CCLV04}
that every unit norm Bessel sequence is a finite union of linearly
independent sets, our next main result gives 
that all unit norm Bessel sequences can be decomposed
into a finite number of sets having Krein's $\omega$-independence.

\begin{theorem}
\label{main2}
Every unit norm Bessel sequence which is finitely linearly independent
is a union of two sets each of which is $\omega$-independent for
$\ell_2$-sequences.  
\end{theorem}

As a main ingredient for the proofs of Theorems \ref{main1} and \ref{main2}, we will
prove a decomposition theorem for frames, which is interesting in
its own right. By providing an explicit construction, we will show that
each unit norm Bessel  sequence can be decomposed into {\em two} subsequences
in such a way that both are small perturbations of ``ideal'' sequences. Our
idea of an ideal sequence is a sequence for which there exists a partition of
its elements into finite sets such that the spans of the elements of
those sets are mutually orthogonal; thus, properties of the sequence are completely
determined by properties of its local components. This definition is inspired by a
more general notion called {\em fusion frames} \cite{CK,CKL06}, which were
designed to model distributed processing applications.


\medskip

This paper is organized as follows. In Section \ref{sec_dabr} we will give the
definition of $\epsilon$-perturbation,
formalize the notion of an ideal sequence, and state some basic results. Section
\ref{sec_decom} contains the Decomposition Theorem and a discussion concerning
an improvement of its proof and concerning the necessity of decomposing into {\it two}
subsequences, whereas the proofs of Theorems \ref{main1} and \ref{main2} will be given in
Section \ref{sec_proof}.


\section{Definitions and basic results}
\label{sec_dabr}

For the remainder let $\H$ be a separable Hilbert space and let $I$ be
a countable index set. Further, in the following we will write $\spa \{f_i\}_{i \in I}$ to mean
the closed linear span of a set of vectors $\{f_i\}_{i \in I}$.

There exist many different definitions for a sequence being a perturbation
of a given sequence. In this paper we will use the following.

\begin{definition}
Let $\{f_i\}_{i \in I}$ and $\{g_i\}_{i \in I}$ be
sequences in $\cH$ satisfying $\{g_i\}_{i \in I} \subset \spa_{i \in I}\{f_i\}$,
and let $\epsilon > 0$. If
\[ \sum_{i \in I}\norm{f_i-g_i}^2 \le \epsilon,\]
then $\{g_i\}_{i \in I}$ is called an {\em $\epsilon$-perturbation of
$\{f_i\}_{i \in I}$}.
\end{definition}

For finite frames, it is precisely the interaction of the frame vectors which makes
them interesting and applicable to a broad spectrum of applications. To employ
the various results which were already obtained in this setting, an ideal infinite
frame should allow its global properties to be determined
locally by considering finite frame sequences.  We formalize this idea in the
next definition in the more general setting of an arbitrary sequence.

\begin{definition}
Let $\{f_i\}_{i \in I}$ be a sequence in $\cH$. We say $\{f_i\}_{i \in I}$
possesses
a {\em finite orthogonal decomposition,} if $I$ can be partitioned
into finite sets
$\{I_j\}_{j = 1}^\infty$ so that
\[ \sp_{i \in I}\{f_i\} =
\left(\sum _{j=1}^\infty \oplus \:  \sp_{i \in I_j}\{f_i\}\right)_{\ell_2}.\]
\end{definition}

We wish to remark that the orthogonal family of finite dimensional subspaces forms
an orthonormal basis of subspaces  and in this sense is a special
case of a {\em Parseval fusion frame} \cite{CK,CKL06}.

The following lemma is well-known, but since the proof is short
and it is fundamental to our construction, we include
it for completeness.

\begin{lemma}
\label{proj_finite}
Let $\{f_i\}_{i \in I}$ be a Bessel sequence in $\cH$, and let $P$ be a finite
rank projection on $\cH$. Then
\[ \sum_{i \in I}\norm{Pf_i}^2 < \infty.\]
\end{lemma}

\begin{proof}
Let $\cK$ be the projection space of $P$ with dimension $d$, and
let $S$ denote the frame operator of $\{Pf_i\}_{i \in I}$, i.e.,
$S(g) = \sum_{i \in I} \ip{g}{f_i}f_i$ for all $g \in \cH$.
Further, let $\{e_j\}_{j=1}^d$ be an orthonormal eigenvector basis
for $\cK$ with respect to $S$ and respective eigenvalues $\{\lambda_j\}_{j=1}^d$.
Then we obtain
\[ \sum_{i \in I}\norm{Pf_i}^2 = \sum_{i \in I}\sum_{j=1}^d\absip{f_i}{e_j}^2
\sum_{j=1}^d\sum_{i \in I}\absip{f_i}{e_j}^2
= \sum_{j=1}^d \lambda_j < \infty.\]
\end{proof}


\section{The Decomposition Theorem}
\label{sec_decom}

The following theorem states that we can decompose each unit norm Bessel sequence into two
subsequences such that both are $\epsilon$-perturbations of sequences which
possess a finite orthogonal decomposition. In fact, we will even derive an explicit
algorithm for generating this partition. The proof is inspired by {\it blocking arguments}
from Banach space theory \cite{LT77}. The Decomposition Theorem will also be the main ingredient
for the proofs of Theorems \ref{main1} and \ref{main2} in Section \ref{sec_proof}.

\begin{theorem}[Decomposition Theorem]
\label{decomp_theo}
Let $\{f_i\}_{i \in I}$ be a unit norm Bessel sequence in $\cH$, and let $\epsilon > 0$.
Then there exists a partition $I = I_1 \cup I_2$ such that, for $j=1,2$, the sequence
$\{f_i\}_{i \in I_j}$ is an $\epsilon$-perturbation of some sequence
in $\cH$, which possesses a finite orthogonal decomposition.
\end{theorem}

\begin{proof}
Let $\{f_i\}_{i \in I}$ be a unit norm Bessel sequence in $\cH$, and let $\epsilon > 0$.
Without loss of generality we may assume that $I = \NN$, since if $I$
is finite we are done.

In the first step we will define a strictly increasing sequence $\{n_i\}_{i  = 1}^\infty$
in $\NN$ by an induction argument.  In the second step, we show that by defining
$I_j := \bigcup_{i = 0}^\infty \{n_{2i+(j-1)}+1,\ldots,n_{2i+j}\}$, for each $j=1,2$,
the sequence $\{f_i\}_{i \in I_j}$ is an $\epsilon$-perturbation of
some sequence $\{g_i\}_{i \in I_j}$ in $\cH$, which possesses a
finite orthogonal decomposition.

For the initial induction step, we set $n_1 := 1$. Further, we define
$S_1$ by $S_1 := \{n_{1}\}$, and
let $P_1$ denote the orthogonal projection onto $\sp_{i \in S_1}\{f_i\}$.
To construct $n_2$, observe that, by Lemma \ref{proj_finite}, we have
\[ \sum_{i = 1}^\infty\norm{P_{1}f_i}^2 < \infty.\]
Therefore we can choose $n_2 > n_1$ so that
\begin{equation}
\label{def_n1}
\sum_{i = n_2+1}^\infty \norm{P_{1}f_i}^2 < \frac{\epsilon}{2}.
\end{equation}
Using this new element of our sequence, we define $T_1$ by $T_1 := \{n_{1}+1,\ldots,n_{2}\}$ and
let $Q_{1}$ denote the orthogonal projection onto $\sp_{i \in T_1}\{f_i\}$.

We proceed by induction. Notice that in each induction step we will define two
new elements of our sequence. Let $k \in \NN$ and suppose that we have already constructed
$n_1,\ldots,n_{2k}$ and defined $\{S_m\}_{m=1}^{k}$, $\{T_m\}_{m=1}^{k}$,
$\{P_m\}_{m=1}^{k}$, and $\{Q_m\}_{m=1}^{k}$. In the following induction step
we will construct $n_{2k+1}$ and $n_{2k+2}$, and define $S_{k+1}$, $T_{k+1}$,
$P_{k+1}$, and $Q_{k+1}$. First, we employ Lemma \ref{proj_finite}, which implies that
\[ \sum_{i = 1}^\infty \norm{Q_{k}f_i}^2 < \infty.\]
Therefore we can choose $n_{2k+1} > n_{2k}$ so that
\begin{equation}
\label{def_n2}
\sum_{i = n_{2k+1}+1}^\infty \norm{Q_{k}f_i}^2 < \frac{\epsilon}{2^{2k}}.
\end{equation}
Now let $S_{k+1}$ be defined by $S_{k+1} := \{n_{2k}+1,\ldots,n_{2k+1}\}$, and let $P_{k+1}$
denote the orthogonal projection of $\cH$ onto $\sp_{i \in \bigcup_{m=1}^{k+1} S_m}\{f_i\}$.
Secondly, again by Lemma \ref{proj_finite}, we have
\[\sum_{i= 1}^\infty \norm{P_{k+1}f_i}^2 < \infty.\]
Thus there exists $n_{2k+2} > n_{2k+1}$ such that
\begin{equation}
\label{def_n3}
\sum_{i = n_{2k+2}+1}^\infty \norm{P_{k+1}f_i}^2 < \frac{\epsilon}{2^{2k+1}}.
\end{equation}
Hence we define the set $T_{k+1}$ by
$T_{k+1} := \{n_{2k+1}+1,\ldots,n_{2k+2}\}$, and let $Q_{k+1}$ denote the orthogonal projection
of $\cH$ onto $\sp_{i \in \bigcup_{m=1}^{k+1} T_m}\{f_i\}$.
Iterating this procedure yields a sequence $\{n_i\}_{i = 1}^\infty$ and, in particular,
we obtain a partition $\{S_m\}_{m= 1}^\infty \cup \{T_m\}_{m= 1}^\infty$
of $\NN$.

For the second step let $\{S_m\}_{m= 1}^\infty$ and $\{T_m\}_{m= 1}^\infty$ be defined as
in the induction argument. Then we define $I_1$ and $I_2$ by
\[I_1 := \bigcup_{m = 1}^\infty S_m \quad \mbox{and} \quad
I_2 := \bigcup_{m = 1}^\infty T_m.\]
It remains to prove that, for $j=1,2$, we can construct a sequence
$\{g_i\}_{i \in I_j}$ in $\cH$ such that $\{g_i\}_{i \in I_j}$ has a
finite orthogonal decomposition and $\{f_i\}_{i \in I_j}$ is an $\epsilon$-perturbation
of it.

In the following we will prove the claim only for $j=1$. The case $j=2$ can be dealt with
in a similar manner.
Using the sequence $\{P_m\}_{m= 1}^\infty$ from the induction argument, we define
$\{g_i\}_{i \in I_1}$ by
\[ g_i := \left\{\begin{array}{rcl}
f_i & : & i \in S_1,\\
f_i - P_{m-1}f_i & : & i \in S_m,\: m > 1.
\end{array}\right.\]
Since by construction we have
\[g_i \in \left\{\begin{array}{rcl}
P_1 \cH & : &i \in S_1,\\
(P_m-P_{m-1})\cH& : & i \in S_m,\: m>1,
\end{array}\right.\]
hence
\[ \sp_{1 \le i < \infty}\{g_i\} =
\left(\sum _{m=1}^\infty \oplus \:  \sp_{i \in S_m}\{g_i\}\right)_{\ell_2},\]
it follows that $\{g_i\}_{i \in I_1}$ possesses a finite orthogonal decomposition.
Further,  for all $m \in \NN$ and $i \in S_m$, we have
\[ P_{m-1}f_i \in \sp_{k \in \bigcup_{l=1}^{m-1} S_l}\{f_k\},\]
which implies that $\sp_{i \in I_1}\{g_i\} = \sp_{i \in I_1}\{f_i\}$.
Finally, applying (\ref{def_n1}) and (\ref{def_n3}) yields
\[ \sum_{i \in I_1} \norm{g_i-f_i}^2 = \sum_{m = 1}^\infty \sum_{i \in S_m}
\norm{g_i-f_i}^2 = \sum_{m >1} \sum_{i \in S_m} \norm{P_{m-1}f_i}^2
\le \sum_{m = 2}^\infty \frac{\epsilon}{2^{m}} < \epsilon.\]
Thus $\{f_i\}_{i \in I_1}$ is an $\epsilon$-perturbation of $\{g_i\}_{i \in I_1}$.
\end{proof}

\begin{remark}  The decomposition argument can be done simultaneously
on two frames at once --- for example on a frame $\{f_i\}_{i \in I}$ and its dual frame,
which is $\{S^{-1}f_i\}_{i \in I}$, $S$ being the frame operator of $\{f_i\}_{i \in I}$.
We will not address this here, since we do not have any serious application
at this time.
\end{remark}

Next we observe that it is necessary to divide our index set into two
subsets in Theorem \ref{decomp_theo}.  That is, the Bessel sequence itself need
not be an $\epsilon$-perturbation of any sequence with a finite
orthogonal decomposition.

\begin{example}\label{E1}
The unit norm Bessel sequence $\{f_i\}_{i = 1}^\infty$ defined by
$f_i = \frac{e_i +e_{i+1}}{\sqrt{2}}$ is not an $\epsilon$-perturbation
of any sequence with a finite orthogonal decomposition for small $\epsilon >0$.
\end{example}

\begin{proof}
If we partition $\NN$ into finite sets $\{I_j\}_{j=1}^{\infty}$, then
there exists a natural number $i_0\in \NN$ so that $i_0\in I_j$ and
$i_0+1 \in I_k$ where $j\not= k$.  Assume, by way of contradiction,
that $\{g_i\}_{i=1}^{\infty}$ is an $\epsilon$-perturbation of
$\{f_i\}_{i=1}^{\infty}$ and $\{g_i\}_{i=1}^{\infty}$ has a finite
orthogonal decomposition given by  $\{I_j\}_{j=1}^{\infty}$.
Then $\norm{f_{i_0}}=1$ and $\norm{f_{i_0}-g_{i_0}} < \sqrt{\epsilon}$, which implies
\[\|g_{i_0} \| \ge \|f_{i_0} \| - \|f_{i_0} - g_{i_0} \| \ge 1-\sqrt{\epsilon}.\]
Similarly, $\|g_{{i_0}+1}\|\ge 1-\sqrt{\epsilon} $.
Since $\spa_{i \in I_j}\{g_i\}$ is orthogonal to $\spa_{i \in I_k}\{g_i\}$, we have
\[\|g_{i_0} - g_{{i_0}+1}\|^2 = \|g_{i_0} \|^2 + \|g_{{i_0}+1}\|^2 \ge 2(1-\sqrt{\epsilon})^2.\]
Using this estimate and the fact that $\|f_{i_0} - f_{{i_0}+1}\|^2 = 1$
and $\|f_{i_0} -g_{i_0} \| + \|f_{{i_0}+1} - g_{{i_0}+1}\| < 2 \sqrt{\epsilon}$,
it follows that
\begin{eqnarray*}
\sqrt{2}(1-\sqrt{\epsilon}) & \le & \|g_{i_0} - g_{{i_0}+1}\|\\
& \le & \|f_{i_0} -f_{{i_0}+1}\| + \|f_{i_0} -g_{i_0} \| + \|f_{{i_0}+1} - g_{{i_0}+1}\|\\
& \le & 1+ 2\sqrt{\epsilon}.
\end{eqnarray*}
This is a contradiction for small $\epsilon >0$.
\end{proof}

\begin{remark}
Our proof of the Decomposition Theorem relies on the ordering of the
elements of the sequence.  This can sometimes cause problems as we
will see below.  However, it is possible to do an
{\it optimal} construction which removes this assumption.  We first
choose $i_{0}\in I$ and let $S_1 = \{i_{0}\}$.  Now, following the
proof,
\[\sum_{i\in I \backslash S_1}\|P_{1}f_i \|^2 < \infty.\]
So choose $T_1 \subset I\backslash S_1$ with $|T_1 |$ minimal and
\[\sum_{i\in I \backslash (T_1 \cup S_1 )} \|P_1 f_i \|^2 < \frac{\epsilon}{2}.\]
So we have put the $f_i$, $i\in I\backslash S_1$ with $\|P_1 f_i \|$ maximal
into $T_1$.  In the induction step (equation (\ref{def_n2})), choose
\[S_{k+1}\subset I \:\backslash \left( \bigcup_{m=1}^{k}S_m \cup \bigcup_{m=1}^{k}T_m
\right)\]
with $|S_{k+1}|$ minimal and
\[\sum_{i \in I \backslash \left( \bigcup_{m=1}^{k+1}S_m \cup \bigcup_{m=1}^{k}T_m
\right)} \|Q_{k} f_i \|^2 < \frac{\epsilon}{2^{2k}}.\]
Similarly, we now construct the next $T_{k+1}$ and then iterate the procedure.
\end{remark}

This stronger form of the decomposition construction is useful because
it eliminates the ordering of the elements.  For example, if we work
with the $\{f_i\}_{i=1}^{\infty}$ in Example \ref{E1}, then for any
permutation of $\{f_i\}_{i=1}^{\infty}$, as long as $f_{i_0} = f_1$
the decomposition we obtain from this stronger form of
the proof of the Decomposition Theorem is
\[\left\{\frac{e_{2i-1}+e_{2i}}{\sqrt{2}}\right\}_{i=1}^{\infty}
\quad \mbox{and} \quad
\left\{\frac{e_{2i}+e_{2i+1}}{\sqrt{2}}\right\}_{i=1}^{\infty}\]
both of which are orthonormal bases for their spans.  But,
if we reorder the sequence $\{f_i\}_{i=1}^{\infty}$ by taking $\{f_i\}_{i=2^k}^{2^{k+1}-1}$ into
\[\{f_{2^k +1}, f_{2^k +2}, \ldots , f_{2^{k+1}-1}, f_{2^k}\} \quad \mbox{for each }
0 \le k < \infty, \]
the proof of the Decomposition Theorem produces the partition
\[I_1 = \{2^{2k},2^{2k+1}, \ldots ,2^{2k+1}-1 : 0 \le k < \infty\},\]
and
\[I_{2} =\{2^{2k+1},2^{2k+1}+1, \ldots ,2^{2k+2}-1 : 0 \le k < \infty\}.\]
Now, $\{f_i\}_{i\in I_j}$ is not even a frame sequence for $j=1,2$.
To see this for $I_1$, we note that the sets
\[J_k = \{2^{2k},2^{2k+1}, \ldots ,2^{2k+1}-1\} \quad (0 \le k < \infty)\]
give a finite orthogonal decomposition of $\{f_i\}_{i\in I_1}$
into linearly independent sets.  So if $\{f_i\}_{i\in I_1}$ would be a
frame sequence, then, since $\{f_i \}_{i\in I_1}$ is ${\omega}$-independent
if and only if for each $k \in \NN$ the sequence $\{f_i \}_{i\in J_k}$ is linearly independent
and by \cite[Proposition 4.3]{Cas00}, it would also be a  Riesz basic sequence.  Let
\[a_{2^{2k}+i} = \frac{(-1)^i}{\sqrt{2^{2k}}},\ \ i=0,1,\ldots , 2^{2k}-1.\]
Then,
\[\sum_{i=0}^{2^{2k}-1}|a_i |^2 = 1,\]
while
\[\left\|\sum_{i=0}^{2^{2k}-1}a_{2^{2k}+i}f_{2^{2k}+i}\right\|^2 =
\left\|\frac{1}{\sqrt{2^{2k}}}(e_1 - e_{2^{2k+1}-1})\right\|^2 = \frac{1}
{2^{2k-1}}.\]
This implies that $\{f_i\}_{i\in I_1}$ is not a Riesz basic sequence. Thus
$\{f_i\}_{i\in I_1}$ is not a frame sequence.


\section{Proofs of Theorems \ref{main1} and \ref{main2}}
\label{sec_proof}

First we will prove Theorem \ref{main2}. For this, we require a small change in the
construction of Theorem \ref{decomp_theo}, which can also be regarded as a strengthening
due to the stronger conditions the decomposition has to satisfy. However, notice that
this modification can only be made provided we have a finitely linearly independent sequence.

\begin{remark}
\label{altering}
We will alter the construction in the proof of Theorem \ref{decomp_theo} for finitely linearly
independent sequences in the following way:

In the $k^{th}$-step, instead of choosing  $n_{2k+1}\ge n_{2k}$
such that
$$
\sum_{i=n_{2k+1}+1}^{\infty}\|Q_{k}f_i\|^2 < \frac{\epsilon}{2^{2k}},
$$
we choose $n_{2k+1}\ge n_{2k}$ so that
$$
\sum_{i=n_{2k+1}+1}^{\infty}\|Q_{k}f_i\|^2 < \frac{\epsilon}{2^{2k}}\delta_k,
$$
where $\delta_k$ denotes the lower Riesz basis bound of $\{f_i\}_{i=1}^{n_{2k}}$.
The choice of $n_{2k+2}\ge n_{2k+1}$ in the same step by \eqref{def_n3} will be adapted
similarly.
\end{remark}

This now enables us to prove Theorem \ref{main2}.

\begin{proof}[{\bf Proof of Theorem \ref{main2}}]
Let $\{f_i\}_{i \in I}$ be a unit norm Bessel sequence which is finitely linearly independent.  
Further, let $\{T_m\}_{m=1}^\infty$, $\{Q_k\}_{k=1}^\infty$, and $\{n_k\}_{k=1}^\infty$ be
chosen as in the proof of Theorem \ref{decomp_theo} modified by Remark \ref{altering}.
We assume that there exists an $\ell_2$-sequence of scalars
$\{a_i\}_{i\in T_m,m=1}^{\ \ \ \ \ \ \infty}$ such that
$$
\sum_{m=1}^{\infty}\sum_{i\in T_m}a_i f_i = 0.
$$
Fix $k\in \N$ and define $g_k$ by
$$
g_k := \sum_{m=1}^{k}\sum_{i\in T_{m}}a_if_i =
- \sum_{m=k+1}^{\infty}\sum_{i\in T_{m}}a_if_i.
$$
Recalling that $Q_{k}$ denotes the orthogonal projection of $\cH$ onto\\
$\sp_{i \in \bigcup_{m=1}^{k} T_m}\{f_i\}$, by the definition of
$g_k$, we have $Q_k g_k = g_k$. Employing this fact, we compute
\begin{eqnarray*}
\|g_k\|^2 &=& - \sum_{m=k+1}^{\infty}\sum_{i\in T_{m}}a_i \langle g_k,f_i\rangle\\
&\le& \left ( \sum_{m=k+1}^{\infty}\sum_{i\in T_{m}} |a_i|^2
\right )^{1/2}\left ( \sum_{m=k+1}^{\infty}\sum_{i\in T_{m}}
|\langle g_k,f_i\rangle|^2\right )^{1/2}\\
&\le& \left ( \sum_{m=k+1}^{\infty}\sum_{i\in T_{m}} |a_i|^2
\right )^{1/2} \|g_k\| \left ( \sum_{m=k+1}^{\infty}\sum_{i\in T_{m}}
\|Q_{k}f_i\|^2 \right )^{1/2}.
\end{eqnarray*}
This implies
\begin{equation} \label{eq:gk}
\|g_k\| \le \left ( \sum_{m=k+1}^{\infty}\sum_{i\in T_{m}} |a_i|^2
\right )^{1/2} \left ( \sum_{m=k+1}^{\infty}\sum_{i\in T_{m}}
\|Q_{k}f_i\|^2 \right )^{1/2}.
\end{equation}
Now let $\delta_k$ denote the lower Riesz basis bound of $\{f_i\}_{i=1}^{n_{2k}}$.
Employing \eqref{eq:gk} and Remark \ref{altering}, we obtain
\begin{eqnarray*}
\delta_k \sum_{m=1}^{k}\sum_{i\in T_{m}}|a_i|^2 &\le&
\|g_k\|^2\\
&\le&
\left ( \sum_{m=k+1}^{\infty}\sum_{i\in T_{m}} |a_i|^2
\right ) \left ( \sum_{m=k+1}^{\infty}\sum_{i\in T_{m}}
\|Q_{k}f_i\|^2 \right )\\
&\le& \left ( \sum_{m=k+1}^{\infty}\sum_{i\in T_{m}} |a_i|^2
\right ) \frac{\epsilon}{2^{2k}}\delta_k.
\end{eqnarray*}
Thus,
$$
\sum_{m=1}^{k}\sum_{i\in T_{m}}|a_k|^2 \le
\frac{\epsilon}{2^{2k}}\sum_{m=k+1}^{\infty}\sum_{i\in T_{m}} |a_i|^2
$$
Since $\{a_i\}_{i\in T_m,m=1}^{\ \ \ \ \ \ \infty}$, 
the left-hand-side of our inequality converges to the $\ell_2$ norm
of this sequence of scalars while the 
right-hand-side converges to zero.  It follows that $a_i=0$ for all
$i=1,2,\ldots$.
\end{proof}

Theorem \ref{main2} will now serve as the main ingredient in the proof of
Theorem \ref{main1}.

\begin{proof}[{\bf Proof of Theorem \ref{main1}}]
Obviously, Conjecture \ref{F_conj} implies Conjecture \ref{our_conj}.

To prove the converse implication suppose that Conjecture \ref{our_conj} holds.
If $\{f_i\}_{i\in I}$ is a unit norm Bessel sequence, it is in particular 
a finite
union of linearly independent sets \cite{CCLV04}. 
Therefore, by Theorem \ref{main2}, it is
a finite union of sets which are $\omega$-independent for
$\ell_2$-sequences.  If $\{f_i\}_{i\in J}$ is one of these families,
by our assumption it is a finite union of frame sequences.  But
each of these frame sequences is $\omega$-independent for $\ell_2$-sequences,
hence it is a Riesz basic sequence.
\end{proof}


\section*{Acknowledgments}
The majority of the research for this paper was performed
while the second and third author were visiting the Department of Mathematics at
the University of Missouri.
These authors thank this department for its hospitality and support
during this visit.
We also thank the American Institute of Mathematics for hosting the workshop
``The Kadison-Singer Problem'' which gave us additional inspirations and also the
possibility to work out the last details of the results presented in this paper
together.


\bibliographystyle{amsplain}

\end{document}